\documentclass{amsart}
\usepackage[latin1]{inputenc}
\usepackage{graphicx}
\usepackage{amssymb}

\usepackage[usenames,dvipsnames]{pstricks}
\usepackage{epsfig}
\usepackage{pst-grad} 
\usepackage{pst-plot} 

\newtheorem{thm}{Theorem}
\newtheorem{cor}[thm]{Corollary}
\newtheorem{lem}[thm]{Lemma}
\newtheorem{prop}[thm]{Proposition}
\theoremstyle{definition}
\newtheorem{defn}[thm]{Definition}
\theoremstyle{remark}

\numberwithin{equation}{section}

\newcommand{\To}{\longrightarrow}

\begin{document}
\setcounter{tocdepth}{1}


\title[]{Finite basis for analytic strong $n$-gaps}
\author{Antonio Avil\'es and Stevo Todorcevic}
\address{Universidad de Murcia, Departamento de Matem\'{a}ticas, Campus de Espinardo 30100 Murcia, Spain.} \email{avileslo@um.es}
\address{CNRS FRE 3233, Universit\'{e} Paris Diderot Paris 7, UFR de math\'{e}matiques case 7012, site Chevaleret
75205 Paris, France. Department of Mathematics, University of Toronto, Toronto, Canada, M5S 3G3.}%
\email{stevo@logique.jussieu.fr, stevo@math.toronto.edu}
\thanks{A. Avil\'{e}s was supported by MEC and FEDER (Project MTM2008-05396), Fundaci\'{o}n S\'{e}neca
(Project 08848/PI/08), Ram\'{o}n y Cajal contract (RYC-2008-02051) and an FP7-PEOPLE-ERG-2008 action.}

\begin{abstract}
We identify the finite list of minimal analytic $n$-gaps which are not weakly countably separated, and we prove that every analytic $n$-gap which is not countably separated contains a gap from our finite list. 
\end{abstract}

\maketitle

\section{Introduction}

We introduced the notion of multiple gap in our recent paper \cite{multiplegaps}. If $\{I_i : i\in n\}$ are $n$-many ideals of $\mathcal{P}(N)$, we say that they constitute a multiple gap (or an $n$-gap) if they are mutually orthogonal\footnote{Notations are explained in Section~\ref{sectionpreliminaries}} and $\bigcap_{i\in n}a_i \neq^\ast\emptyset$ whenever $a_i\geq I_i$ for every $i$. Such a multiple gap is weakly countably separated if there exists a countable family $\mathcal{C}$ of infinite subsets of $N$ such that whenever we choose $x_i\in I_i$ we can find $\{c_i : i\in n\}\subset \mathcal{C}$ such that $\bigcap_{i\in n}c_i  =^\ast \emptyset$ and $x_i\subset^\ast c_i$ for each $i\in n$.\\

In this paper we are interested in multiple gaps which are not weakly countably separated, so for short we shall call them \emph{strong} multiple gaps. When we denote an $n$-gap by a letter $\Gamma$, we understand that $\Gamma = \{\Gamma(i) : i\in n\}$ and each $\Gamma(i)$ is an ideal. An $n$-gap $\Gamma$ is analytic if each of the ideals $\Gamma(i)$ is an analytic family of subsets of $N$ (when we view $\Gamma(i)\subset 2^N$ and identify $2^N$ with the Cantor set). We define an order relation between multiple gaps:\\

\begin{defn}
Suppose that $\Gamma$ and $\Gamma'$ are $n$-gaps on the countable sets $N$ and $N'$ respectively. We say that $\Gamma\leq \Gamma'$ if there is a one-to-one map $\phi:N\To N'$ such that for every $x\in\mathcal P(N)$,
\begin{enumerate}
\item $\phi(x)\in \Gamma'(i)$ whenever $x\in \Gamma(i)$, $i\in n$.
\item $\phi(x) \in (\Gamma')^\perp$ whenever $x\in \Gamma^\perp$.\\
\end{enumerate}
\end{defn}

\begin{defn}
A strong analytic $n$-gap $\Gamma$ is called minimal if for every analytic strong $n$-gap $\Gamma'$, if $\Gamma'\leq \Gamma$ then $\Gamma\leq \Gamma'$. Two minimal strong $n$-gaps $\Gamma$ and $\Gamma'$ are called equivalent if $\Gamma\leq\Gamma'$, equivalently $\Gamma'\leq\Gamma$.\\
\end{defn}

We obtain the following result:

\begin{thm}\label{existenceofminimals}
For a fixed finite number $n$, there exists only a finite number of equivalence classes of analytic minimal strong $n$-gaps. Moreover, for every analytic strong $n$-gap $\Gamma$, there exists a minimal one $\Gamma'$ such that $\Gamma'\leq \Gamma$.
\end{thm}

This is proven by combining, on the one hand our multidimensional generalization \cite[Theorem 7]{multiplegaps} of the second dichotomy for analytic gaps from \cite{Todorcevicgap}, and on the other hand Milliken's Ramsey theorem for trees \cite{Milliken}. The concrete form of the list of minimal analytic strong $n$-gaps is provided by Theorem~\ref{thetypes} in Section~\ref{sectiontypes}. This list is not obtained just as a by-product of the proof of Theorem~\ref{existenceofminimals}, but it requires a further minute combinatorial analysis to which most of this paper is devoted. In the cases when $n=2,3$, for example, we get that there exist (up to equivalence and up to permutation) four analytic minimal strong 2-gaps and nine analytic minimal strong 3-gaps which are described at the end of this paper. We shall also give an application of these results, by proving that an analytic strong multiple gap cannot be a clover, giving a partial solution to a problem that we posed in \cite{multiplegaps}.

\section{Preliminaries}\label{sectionpreliminaries}

\textbf{Set operations}. We identify a natural number with the set of its predecessors: $n=\{0,1,\ldots,n-1\}$; $a\subset^\ast b$ means that $b\setminus a$ is finite, and $a=^\ast b$ means that $a\subset^\ast b$ and $b\subset^\ast a$. We say that two sets are orthogonal if $a\cap b =^\ast \emptyset$ and we say that two families $I$ and $J$ of sets are orthogonal if $a\cap b =^\ast \emptyset$ whenever $a\in I$, $b\in J$. If $I$ is a family of subsets of $N$, we write $I^\perp = \{a\subset N : \forall b\in I\ a\cap b =^\ast \emptyset\}$. We write $a\geq I$ if $a\supset^\ast b$ for all $b\in I$. If $\Gamma = \{\Gamma(i): i\in n\}$ is an $n$-gap, we write $\Gamma^\perp = (\bigcup_{i\in n}\Gamma(i))^\perp$. An $n$-gap is dense if $\Gamma^\perp$ consists only of finite sets. We often work with the family $\mathcal{P}(N)$ of all subsets of a suitable countable set $N$, identifying it with $\mathcal{P}(\omega)$. An ideal of $\mathcal{P}(N)$ is a nonempty $I\subset \mathcal{P}(N)$ which is closed under finite unions and under passing to arbitrary subsets. We are only interested in ideals of $\mathcal{P}(N)$ which correspond to ideals of $\mathcal{P}(N)/fin$, so we always assume that an ideal contains all finite sets.\\

\textbf{The $m$-adic tree}. For a natural number $m$, we consider $m^{<\omega}$ the $m$-adic tree consisting of finite sequences of elements of $m$, endowed with the tree order: $s\leq t$ if $s=(s_0,\ldots,s_k)$, $t=(t_0,\ldots,t_l)$, $k\leq l$ and $s_i = t_i$ for $i\leq k$. Given any $s=(s_0,\ldots,s_k)$ and $t=(t_0,\ldots,t_l)$ in $m^{<\omega}$. we denote $s^\frown t = (s_0,\ldots,s_k,t_0,\ldots,t_k)$, while $r = s\wedge t$ is the maximal element such that $r\leq s$ and $r\leq t$. If $i\in m$, then $s^\frown i = s^\frown (i)$. If $s,t\in m^{<\omega}$, and $s^\frown i \leq t$ ($i \in m$), then we call $inc(s,t) = inc(t,s) = (i,i)$. Otherwise, if there is an $r = s\wedge t$ such that $r^\frown i = s$ and $r^\frown j = t$, we call $inc(s,t) = (i,j)$. We also denote $|(s_0,\ldots,s_{k-1})| = k$.\\

\textbf{Chains and combs}. Given $u,v\in m$, a subset $a\subset m^{<\omega}$ is called a $u$-chain if $a = \{s[0],s[1],\ldots\}$ with $s[i]^\frown u \leq s[i+1]$ for every $i$. If $u\neq v$, $u,v\in m$, a set $b= \{t[0],t[1],\ldots\}\subset m^{<\omega}$ is called a $(u,v)$-comb if there exists a $u$-chain $a = \{s[0],s[1],\ldots\}$ such that $s[i]^\frown v \leq t[i]$ and $|t[i]| < |s[i+1]|$ for every $i$. For convenience, we define $(u,u)$-combs to be $u$-chains.\\

\textbf{Ramsey theory of complete subtrees}. A subset $T\subset m^{<\omega}$ is a complete subtree if $T$ is the range of a one-to-one function $\phi:m^{<\omega}\To m^{<\omega}$ such that $s^\frown u \leq t$ implie $\phi(s)^\frown u \leq \phi(t)$ for every $s,t\in m^{<\omega}$ and every $u\in m$. We shall need the following statement:

\begin{thm}\label{Ramseytree}
Let $T$ be a complete subtree of $m^{<\omega}$. Let $u,v\in m$, let $\mathcal{C}$ be the set of all $(u,v)$-combs, $k$ a natural number, and $c:\mathcal{C}\To k$ a Souslin-measurable function. Then there exists a complete subtree $T'\subset T$ such that $c$ takes the same value on all $(u,v)$-combs of $T'$.
\end{thm}

Proof: The application of Milliken's theorem \cite[Theorem 6.13]{Ramsey} provides a complete subtree $T_0\subset T$ whis is moreover \emph{strong subtree} (meaning that that the function $\phi$ above takes elements on the same level to elements on the same level) and such that $c$ takes the same value on all $(u,v)$-combs of $T_0$ of the same \emph{embedding type} \cite[Definition 6.11]{Ramsey}. However, we can get a complete subtree $T_1\subset T_0$ such that all $(u,v)$-combs of $T_1$ have the same embedding type inside $T_0$. For this, just make sure that for every level $l$, there is at most one element $t[l]\in T_1$ in the $l$-th level of $T_0$, and for every $s=(s_0,\ldots,s_k)\in T_1$ at an upper level, we have either $s|_l = t[l]$ or $s_l = 0$.$\qed$\\

\section{The proof of Theorem \ref{existenceofminimals}}

\begin{defn}
Given a function $f:m^2\To n\cup\{\infty\}$, and $i\in n$, we consider the ideal $\Gamma_f(i)$ on $m^{<\omega}$ generated by all $(u,v)$-combs such that $f(u,v) = i$. We also denote $\Gamma_f = \{\Gamma_f(i) : i\in n\}$.
\end{defn}

\begin{prop}
If the range of $f$ contains $n$, then $\Gamma_f = \{\Gamma_f(i) : i\in n\}$ is a strong $n$-gap.
\end{prop}

Proof: First, notice that the ideals $\{\Gamma_f(i) : i\in n\}$ are mutually orthogonal, because the intersection of a $(u,v)$-comb and a $(u',v')$-comb is finite whenever $(u,v)\neq(u',v')$. If we consider the function $g:n^2\To n\cup\{\infty\}$ given by $g(i,i) = i$ and $g(i,j) = \infty$ if $i\neq j$, then $\Gamma_g(i)$ is the ideal generated by $i$-chains. It is proven in \cite[Theorem 6]{multiplegaps} that in such case, $\Gamma_g$ is a strong $n$-gap. We will reduce the case of a general $f$ to this particular one. 
For every $i\in n$ choose $(u_i,v_i)$ such that $f(u_i,v_i)=i$. We consider $s\in m^{<\omega}$ and elements $s[i]\in m^{<\omega}$ for $i\in n$ all of the same length, greater than the length of $s$, such that $inc(s[i],s) = (u_i,v_i)$, as shown in the picture:\\


\begin{center}
\includegraphics[scale=0.5]{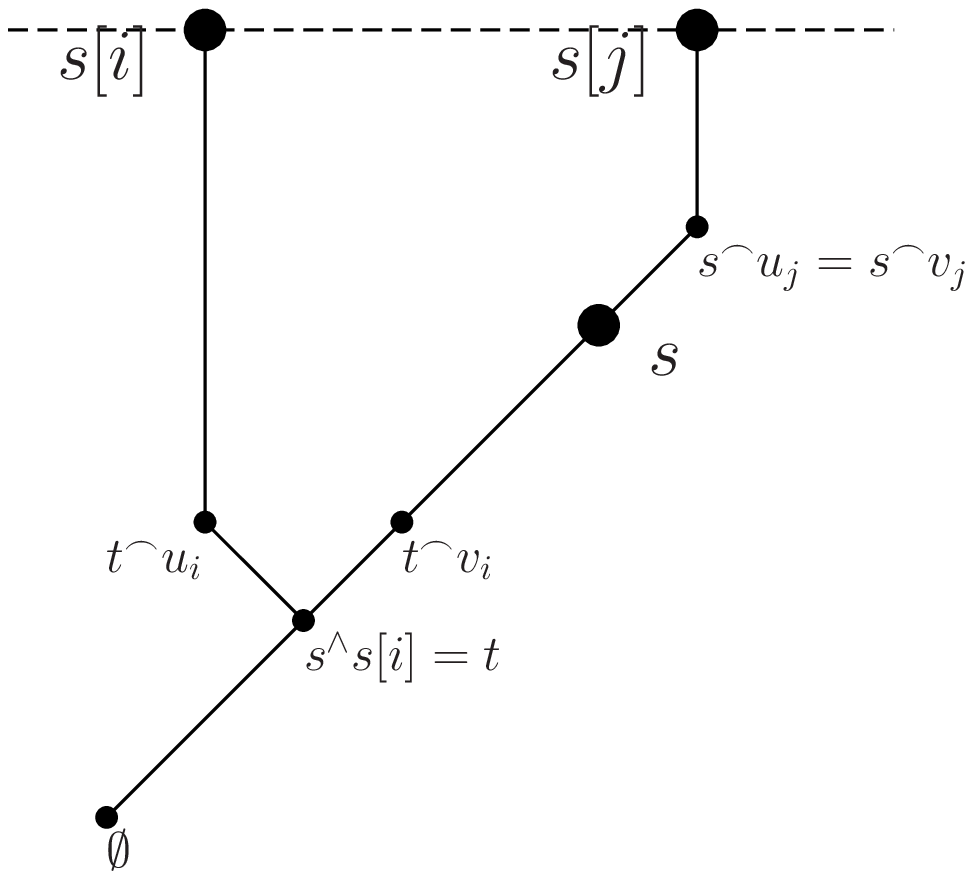}
\end{center}


Define a function $\phi:n^{<\omega}\To m^{<\omega}$ by $$\phi(j_0,j_1,\ldots,j_k) = s[j_0]^\frown s[j_1]^\frown\cdots^\frown s[j_k]^\frown s.$$
This is a one-to-one function such that $\phi(A)$ is a $(u_i,v_i)$-comb whenever $A$ is an $i$-chain. Hence, we have that $\phi(\Gamma_g(i))\subset \Gamma_f(i)$. Since we knew that $\Gamma_g$ is a strong $n$-gap, we get that $\Gamma_f$ is also a strong $n$-gap.$\qed$\\

\begin{lem}\label{densityofcombs}
For every infinite set $x\subset m^{<\omega}$ there exists an infinite subset $x'\subset x$ such that $x'$ is a $(u,v)$-comb for some pair $(u,v)\in m^2$.
\end{lem}

Proof:We construct sequences $s_1,t_1,s_2,t_2,\ldots\in m^{<\omega}$ and $x= x_1\supset x_2\supset\cdots$ such that $|s_i|<|t_i|<|s_j|$, $s_i<s_j\wedge t_i$, $t_i\in x_i$ and $\{t\in x_i : t>s_i\}$ is infinite for all $i<j$. This sequences are defined inductively, namely, let $x_1=x$, $s_1=\emptyset$, and $t_1\in x$ arbitrary; given $x_i$, $s_i$ and $t_i$, fix an integer $k_i>|t_i|$ and find $s_{i+1}>s_i$ with $|s_{i+1}|=k_i$ such that $x_{i+1} = \{t\in x_{i} : t>s_{i+1}\}$ is infinite, and then choose any $t_{i+1} \in x_{i+1}$. Once these sequences are defined, since $s_1<s_2<\cdots$ and $t_i>s_i$, it is clear that there is an infinite set $A$ such that $x'=\{t_i : i\in A\}\subset x$ is a $(u,v)$-comb for some $(u,v)$.$\qed$\\

\begin{lem}\label{density}
For every function $f:m^2\To n\cup\{\infty\}$, $\Gamma_f(\infty)$ is a dense subideal of $\Gamma_f^\perp$. Hence, $\Gamma_f$ is dense if and only if $\infty$ is not in the range of $f$.
\end{lem}

Proof: It is clear that $\Gamma_f(\infty)\subset \Gamma_f^\perp$, again becuase the intersection of a $(u,v)$-comb and a $(u',v')$-comb is finite whenever $(u,v)\neq (u',v')$. Being dense means that for every infinite set $x\in\Gamma_f^\perp$ there exists an infinite $x'\in\Gamma_f(\infty)$ such that $x'\subset x$. If $x'$ is a $(u,v)$-comb, as given by Lemma~\ref{densityofcombs}, then we have $x'\in\Gamma_f(\infty)$ as desired, because $f(u,v)=\infty$ (otherwise, if $f(u,v)=i\neq\infty$, then we would have that $x'\in \Gamma_f(i)$ but $x'\subset x\in \Gamma_f^\perp$).$\qed$\\

\begin{thm}\label{strongcontain}
Let $\Gamma$ be an analytic strong $n$-gap. Then there exists $f:n^2\To n\cup\{\infty\}$ such that $\Gamma_f$ is an $n$-gap and $\Gamma_f\leq \Gamma$.
\end{thm}

Proof: We consider again the function $g:n^2\To n\cup\{\infty\}$ given by $g(i,i) =i$ and $g(i,j) = \infty$ for $i\neq j$. We proved in \cite[Theorem 7]{multiplegaps} that there exist a one-to-one function $\varphi:n^{<\omega}\To N$ such that $\varphi(x)\in \Gamma(i)$ whenever $x\in \Gamma_g(i)$, $i\in n$. Given a pair $(i,j)$ we can consider a function $c$ that associates to each $(i,j)$-comb $x\subset n^{<\omega}$ a $c(x)\in n\cup\{\infty\}$ such that $\varphi(x)\in \Gamma(c(x))$ if $c(x)\neq \infty$, or $\varphi(x)\not\in\bigcup_{k\in n}\Gamma(k)$ if $c(x)=\infty$. By repeated application of Theorem~\ref{Ramseytree}, we can find a complete subtree $T$ of $n^{<\omega}$ such that for every $(i,j)$ the corresponding function $c$ is constant equal to some $h(i,j)$ on the set of all $(i,j)$-combs of $T$. Consider the function $f$ given by $f(i,i) = g(i,i)$ and $f(i,j) = h(i,j)$ for $i\neq j$. Let us see that $\varphi:T\To N$ witnesses that $\Gamma_f\leq \Gamma$. Indeed, if $x\in \Gamma_f(\xi)$ is an $(i,j)$-comb of $T$, so that $f(i,j)=\xi\neq\infty$, then $\varphi(x)\in \Gamma(\xi) = \Gamma(c(x))$. On the other hand, suppose that we had $x\subset T$ such that $x\in\Gamma_f^\perp$ but $\varphi(x)\not\in \Gamma^\perp$. Then, since $\varphi(x)\not\in \Gamma^\perp$, there is an infinite set $y\subset\varphi(x)$ such that $y\in \Gamma(\xi)$ for some $\xi\in n$. Then by Lemma~\ref{density}, since $\varphi^{-1}(y)\subset x\in \Gamma_f^\perp$, there is an infinite set $x'\subset \varphi^{-1}(y)$ such that $x'\in\Gamma_f(\infty)$, and we can suppose that $x'$ is an $(i,j)$-comb with $f(i,j)=\infty$. Then $c(x')=\infty$, which implies that $\varphi(x')\not\in \bigcup_{k\in n}\Gamma(k)$, and this contradicts that $\varphi(x')\subset y \in \Gamma(\xi)$.   $\qed$\\

Theorem \ref{existenceofminimals} is now a corollary of Theorem~\ref{strongcontain}, which provides a finite set of analytic strong $n$-gaps (those of the form $\Gamma_f$) which are below any other. Then minimal elements of this finite set under the relation $(\leq)$ are indeed minimal strong $n$-gaps and are below any other analytic strong $n$-gap. Thus, the minimal analytic strong $n$-gaps are the strong $n$-gaps of the form $\Gamma_f$ which are minimal among the $n$-gaps of the same form. The rest of the paper is devoted to this task. The first step is to understand when $\Gamma_f \leq \Gamma_g$.\\

\section{Analyzing when  $\Gamma_f \leq \Gamma_g$}

\begin{defn}
For natural numbers $m_0$ and $m_1$, we say that $\varepsilon:m_0^2\To m_1^2$ is a reduction map if there exists $k<\omega$ (which can be taken less than $m_0+2$), a one-to-one function $e:m_0 \To m_1^k$ and an element $x\in m_1^{<k}$ such that
\begin{enumerate}
\item $\varepsilon(u,v) = inc(e(u),e(v))$ for $u\neq v$,
\item $\varepsilon(u,u) = inc(e(u),x)$.\\
\end{enumerate}
\end{defn}

\begin{lem}\label{criterion}
Suppose that we have two $n$-gaps of the form $\Gamma_f$ and $\Gamma_g$, for functions $f:m_0^2\To n\cup\{\infty\}$ and $g:m_1^2\To n\cup\{\infty\}$. Then the following are equivalent:
\begin{enumerate}
\item $\Gamma_f\leq\Gamma_g$.
\item There exists a reduction map $\varepsilon:m_0^2\To m_1^2$ such that $f = g\circ \varepsilon$.
\end{enumerate}
\end{lem}

Proof: For $(2)\Rightarrow (1)$, consider the reduction map $\varepsilon$ and its associated $e:m_0 \To m_1^k$ and $x\in m_1^{<k}$. Define a function $\phi: m_0^{<\omega}\To m_1^{<\omega}$ by
$$\phi(u_0,\ldots,u_k) = e(u_0)^\frown e(u_1)^\frown\cdots^\frown e(u_k)^\frown x.$$
Notice that the function $\phi$ maps $(u,v)$-combs onto $\varepsilon(u,v)$-combs. We prove that this property implies that $\phi$ actually witnesses that $\Gamma_f\leq \Gamma_g$. If $x\in\Gamma_f(i)$, we can suppose that $x$ is a $(u,v)$-comb where $f(u,v)=i$, then $\phi(x)$ is an $\varepsilon(u,v)$-comb, where $g(\varepsilon(u,v)) = f(u,v) = i$, hence $\phi(x)\in\Gamma_g(i)$. We check now that if $x\in \Gamma_f^\perp$, then $\phi(x)\in\Gamma_g^\perp$. If $\phi(x)\not\in\Gamma_g^\perp$, then by Lemma~\ref{densityofcombs} there exists $y\subset\phi(x)$ a $(u',v')$-comb, where $g(u',v')=i\neq\infty$. Again by Lemma~\ref{densityofcombs}, we find $x'\subset\phi^{-1}(y)\subset x$ a $(u,v)$-comb. Then $\phi(x')\subset y$ is a $\varepsilon(u,v)$-comb and an $(u',v')$-comb as well, so $(u',v') = \varepsilon(u,v)$. Hence $f(u,v) = g(\varepsilon(u,v)) = g(u',v') = i\neq\infty$ contradicting that $x'\subset x\in \Gamma_f^\perp$.\\

We prove now that $(1)$ implies $(2)$. Suppose $\Gamma_f\leq \Gamma_g$ and let $\phi:m_0^{<\omega}\To m_1^{<\omega}$ be a one-to-one map witnessing this fact. Along the proof, we need to use the fact that if a complete subtree is divided into finitely many pieces, then one of them contains a further complete subtree. This is easy to prove: if $T=P\cup Q$, then either for every $t\in P$ there are elements of $P$ above each $t^\frown i$ (in which case $P$ contains a complete subtree) or else $Q$ contains a complete subtree.\\

We define a complete subtree $T\subset m_0^{<\omega}$. The construction is done inductively, so that given $t\in T$ we define its immediate successors. Also, for every $t\in T$ we define along the construction a complete subtree $T_t$ rooted at $t$. We start by declaring $\emptyset\in T$, and $T_\emptyset = m_0^{<\omega}$. So, given $t\in T$, we choose a level $l$ of the tree $m_1^{<\omega}$ higher than the level of $\phi(t)$, and for every $i\in m_0$ an element\footnote{we can find these $s_i(t)$ and $T_i$ by the above mentioned property that if a complete subtree is divided into finitely many pieces, one of them contains a further complete subtree.}  $s_i(t)$ on the $l$-level of $m_1^{<\omega}$ for which there is a complete subtree $T_{i}\subset T_t$ above $t^\frown i$ so that $s_i(t)\leq \phi(r)$ for every $r\in T_i$.  The level $l$ and the $s_i(t)$ are chosen so that the cardinality of $\{s_i(t) : i\in m_0\}$ is maximized among all possible choices that we could make. The immediate successor of $t$ in $T$ are chosen to be the roots $t[i]$ of every $T_i$, and $T_{t[i]} = T_i$.\\

\begin{center}
\includegraphics[scale=0.75]{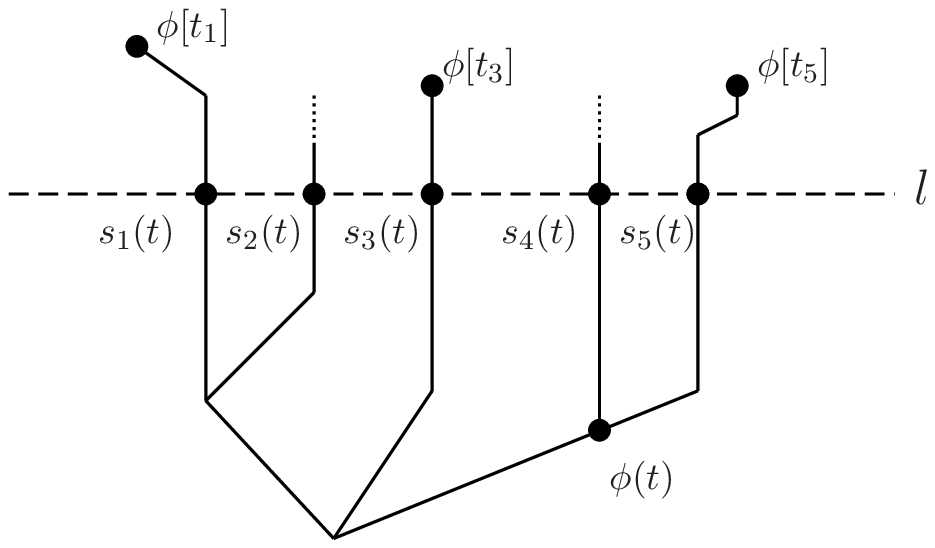}
\end{center}

By passing to a further complete subtree, we can suppose that for each $t\in T$, $\phi(t)$ and the $s_i(t)$ follow the same pattern of incidences (meaning that $inc(\phi(t),s_i(t))$ and $inc(s_i(t),s_j(t))$ are independent of $t$).\\

Claim A: $|\{s_i(t) : i\in m_0\}|$ equals either $m_0$ or 1. Proof: suppose that we have $x = s_i(t) = s_j(t)$ for some $i\neq j$, but also $s_p(t)\neq s_q(t)$ for some other $p,q$. Because the pattern of incidences is independent of $t$, we have that $s_p(t[i])\neq s_q(t[i])$ and $s_p(t[j])\neq s_q(t[j])$. It follows that for some $\xi,\zeta\in\{p,q\}$, $s_\xi(t[i])$ and $s_\zeta(t[j])$ are incomparable nodes. But then, there is a complete subtree $T'_i\subset T_t$ above $t[i]^\frown\xi > t^\frown i$ all of whose nodes are mapped by $\phi$ above $s_\xi(t[i])$, and a complete subtree $T'_j\subset T_t$ above $t[j]^\frown \zeta>t^\frown j$ all of whose nodes are mapped by $\phi$ above $s_\zeta(t[j])$. This implies that if we take $l'$ high enough, we could have chosen nodes of level $l'$ in our construction so that $s'_i(t)>s_\xi(t[i])$, $s'_j(t)>s_\zeta(t[j])$, $s'_k(t)>s_k(t)$ for all other $k$'s, as shown in the picture\footnote{Again, we know that these nodes $s'_k(t)$ can be found by the principle that if a complete subtree is divided into finitely many pieces, one of them contains a complete subtree}:  

\begin{center}
\includegraphics[scale=0.75]{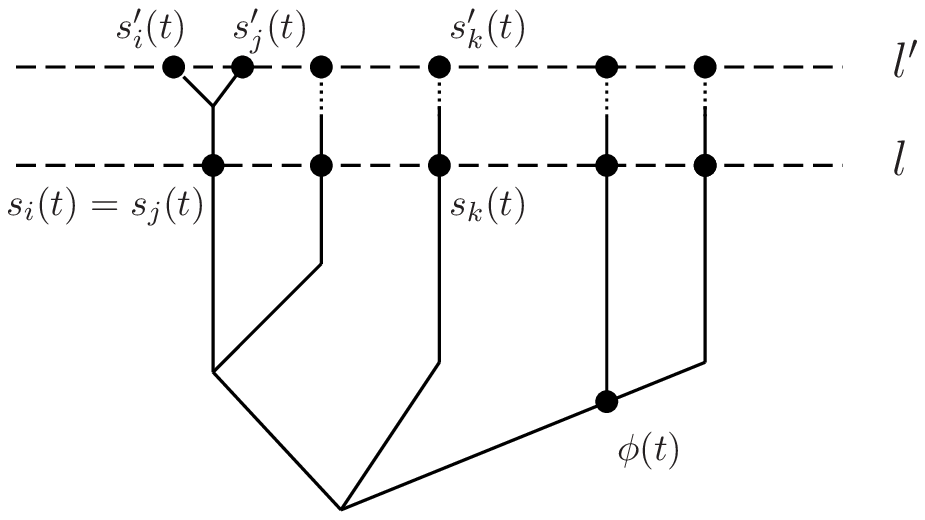}
\end{center}

This means that we could have chosen $l'$ instead of $l$ and $s'_k(t)$ instead of $s_k(t)$ and we would have that the cardinality of $\{s'_k(t): k\in m_0\}$ is one unit larger than that of $\{s_k(t) : k\in m_0\}$. This contradicts the rules of our construction and finishes the proof of Claim A.\\

Claim B: $|\{s_i(t) : i\in m_0\}| = m_0$. If not, by Claim A we would have that $|\{s_i(t) : i\in m_0\}| = 1$ for all $t$, call $s(t) = s_i(t)$. We prove that the set $\{s(t) : t\in T\}$ is contained in a branch $x$ of $m_1^{<\omega}$. Indeed, suppose that $s(t)$ and $s(t')$ are incomparable for $t\neq t'$. It is easy to see that we can even suppose that $t$ and $t'$ are incomparable. Let $r=t\wedge t'$, $t\geq r^\frown u$ and $t'\geq r^\frown v$. Then, we can find a complete subtree $T'\subset T_r$ above $r^\frown u$ all of whose elements are mapped by $\phi$ above $s(t)$, and a complete subtree $T''\subset T_r$ above $r^\frown v$ all of whose elements are mapped by $\phi$ above $s(t')$. In a similar way as we did before, this would imply the possibility to have chosen a larger level $l'$ and different $s'_u(r)> s(t)$ and $s'_v(r)>s(t')$, a contradiction. So we suppose that all $s(t)$ are inside the branch $x$. For every $t\in T$ consider $r_t\in x$ with $|r_t|>|\phi(t)|$ and define $inc(x,\phi(t)) = inc(r_t,\phi(t))$ which is independent of the choice of $r_t$. By passing to a nice subtree $T'\subset T$ we can suppose that $inc(x,\phi(t)) =(u,v)$ is independent of $t$, and we can even suppose that $\phi(T')$ is just a $(u,v)$-comb. But then it is impossible that $\phi$ witnesses $\Gamma_f\leq\Gamma_g$ because it would send any type of combs of $T'$ to the same type of $(u,v)$-comb.\\

Consider the reduction map $\varepsilon$ associated to the incidence pattern of the $s_i(t)$ and $\phi(t)$. That is,
\begin{enumerate}
\item $\varepsilon(i,j) = inc(s_i(t),s_j(t))$ for $i\neq j$,
\item $\varepsilon(i,i) = inc(s_i(t),\phi(t))$.
\end{enumerate}

We observe that $i$-chains of $T$ map under $\phi$ into $\varepsilon(i,i)$-chains or combs, and $(i,j)$-combs map under $\phi$ into $\varepsilon(i,j)$-combs. So now, let $a$ be an $(i,j)$-comb in $m_0^{<\omega}$, then $a\in \Gamma_f(f(i,j))$. On the one hand, since $\phi$ withness that $\Gamma_f\leq \Gamma_g$ we get that $\phi(a)\in \Gamma_g(f(i,j))$. On the other hand, since $\phi(a)$ is an $\varepsilon(i,j)$-comb, we get that $\phi(a)\in \Gamma_g(g(\varepsilon(i,j))$. Therefore, $f(i,j) = g(\varepsilon(i,j))$.$\qed$\\  

\section{Finding invariants}

The next thing that we are going to do is to find some characteristics of a function $f:m^2\to n\cup\{\infty\}$ that will be well behaved under the relation $\Gamma_f\leq \Gamma_g$. These characteristics will be invariant under equivalence of minimal strong $n$-gaps, thus giving us the first hint of how many equivalence classes exist. Let us call $\Delta(m) = \{(i,i) : i\in m\}$.

\begin{defn} Let $f:m^2\To n\cup\{\infty\}$, 
\begin{enumerate}
\item $pbranch(f) = \{k\in n : f^{-1}(k)\subset\Delta(m)\}$.
\item Given $k,l\in n\cup\{\infty\}$, we say that $k$ is $f$-attached to $l$ if: $f(i,j) = k$, $i\neq j$, implies that $f(j,i)=l$.\\
\end{enumerate}
\end{defn}

From this moment on, when we say ``$\Gamma_f$ is an $n$-gap'' it is implicitly understood that $f:m^2\To n\cup\{\infty\}$ is a function whose range covers $n$, for some nautral number $m$.\\ 

\begin{lem}\label{hereditarypbranch}\label{attachpreserved}
Suppose $\Gamma_f\leq \Gamma_g$ are $n$-gaps, then
\begin{enumerate}
\item $pbranch(g) \subset pbranch(f)$.
\item For every $k,l\in n\cup\{\infty\}$, if $k$ is $g$-attached to $l$, then $k$ is $f$-attached to $l$.
\end{enumerate}
\end{lem}

Proof: Part (1) follows from the fact that if $\varepsilon:m_0^2\To m_1^2$ is a reduction map, then the two coordinates of $\varepsilon(u,v)$ are different whenever $u\neq v$. The other item follows from the fact that for every reduction map, if $\varepsilon(u,v) = (i,j)$, then $\varepsilon(v,u) = (j,i)$.$\qed$\\ 

\begin{cor}
If $\Gamma_f$ and $\Gamma_g$ are equivalent minimal $n$-gaps, then 

\begin{enumerate}
\item $pbranch(g) = pbranch(f)$.
\item For every $k,l\in n\cup\{\infty\}$, ($k$ is $g$-attached to $l$) $\Leftrightarrow$ ($k$ is $f$-attached to $l$).\\
\end{enumerate}
\end{cor}

\begin{lem}\label{attachproperties}
Let $k,l,m\in n\cup\{\infty\}$, $k\not\in pbranch(f)$. 
\begin{enumerate}
\item If $k$ is $f$-attached to $l$, and $k$ is $f$-attached to $m$, then $l=m$.
\item If $k$ is $f$-attached to $l$, and $l$ is $f$-attached to $m$, then $k=m$.
\end{enumerate}
\end{lem}

Proof: Since $k\not\in pbranch(f)$, there exist $i\neq j$ with $f(i,j)=k$. In the first case, we have $l = f(j,i) = m$. In the second case we have $f(j,i) = l$, and then $k = f(i,j) = m$.$\qed$.

\begin{lem}\label{pbranchdychotomy}
Let $\Gamma_g$ be an $n$-gap. Then one and only one of the two following possibilities holds.
\begin{enumerate}
\item Either: $pbranch(g)=\emptyset$ and for every $k\in n$ there exists $l\in n$ such that $l\neq k$ and $l$ is $g$-attached to $k$;
\item Or: there exists an $n$-gap $\Gamma_f$ such that $\Gamma_f\leq \Gamma_g$ and $pbranch(f)\neq\emptyset$.
\end{enumerate}
\end{lem}

Proof:  We have $g:m_1^2\To n\cup\{\infty\}$. If $pbranch(g)\neq\emptyset$, then (2) holds, so suppose that $pbranch(g)=\emptyset$ but the first possibility does not hold. We can suppose without loss of generality that there is no $l\in n\setminus\{0\}$ which is $g$-attached to $0$. Combining this with the fact that $pbranch(g)=\emptyset$, we can choose for every $k\in n$ a pair $(i_k,j_k) \in m_1^2\setminus \Delta(m_1)$ such that $g(i_k,j_k) = k$ and $l(k) : = g(j_k,i_k) \neq 0$ if $k>0$. We consider a reduction map $\varepsilon:m_0^2 = n^2\To m_1^2$ corresponding to a function $e:n\To m_1^n$ and an $x$ described as follows:

\begin{eqnarray*}
x &=& (j_0);\\
e(k) &=& (i_0,\ldots,i_k,j_{k+1},0,\ldots,0), k=0,\ldots,n-2;\\
e(n-1) &=& (i_0,\ldots,i_{n-1}).\\
\end{eqnarray*}

Consider then $f(u,v) = g(\varepsilon(u,v))$, so that we will have $\Gamma_f\leq\Gamma_g$.  Then
\begin{itemize}
\item If $k<k'<n$, then $f(k',k) = g(inc(e(k'),e(k)) = g(i_{k+1},j_{k+1}) = k+1$
\item If $k<k'<n$, then $f(k,k') = g(inc(e(k),e(k')) = g(j_{k+1},i_{k+1}) = l(k+1)$
\item If $k<n$, then $f(k,k) = g(inc(e(k),x)) = g(i_0,j_0) = 0$
\end{itemize}
It follows that $\Gamma_f$ is really an $n$-gap, since all elements of $n$ are in the range of $f$, and also $0\in pbranch(f)$.$\qed$\\

In the sequel, we use the notation $\langle A \rangle^2 = \{(i,j)\in A^2 : i\neq j\}$.\\

\begin{defn}
Let $\Gamma_f$ be an $n$-gap, let $A\subset pbranch(f)$, and $\psi:\langle A \rangle^2\To n\cup\{\infty\}\setminus A$. We call $f$ to be $\psi$-branch-reduced if $f(i,j) = \psi(f(i,i),f(j,j))$ whenever $f(i,i),f(j,j)\in A$ and $f(i,i)\neq f(j,j)$.\\
\end{defn}

\begin{lem}\label{reducedpreserved}
Let $\Gamma_g$ be an $n$-gap, let $A\subset pbranch(g)$ and $\psi:\langle A \rangle^2\To n\cup\{\infty\}\setminus A$. If $\Gamma_f\leq \Gamma_g$ and $\Gamma_g$ is $\psi$-branch-reduced, then $f$ is also $\psi$-branch-reduced.\\
\end{lem}

Proof: We know that when $g(i,i),g(j,j)\in A$, $g(i,i)\neq g(j,j)$, then $g(i,j) = \psi(g(i,i),g(j,j))$. Suppose that we are given $u,v$ such that $f(u,u),f(v,v)\in A$ and they are different. Consider a reduction map $\varepsilon$ that gives $\Gamma_f\leq \Gamma_g$, with its associated function $e$ and element $x$. Since $f(u,u) = g(\varepsilon(u,u)) \in A\subset pbranch(g)$, there exists $i$ such that $\varepsilon(u,u) = (i,i)$ and similarly, there exists $j$ with $\varepsilon(v,v) = (j,j)$; (notice also that $(i,i)\neq (j,j)$ since $g(i,i) = f(u,u)\neq f(v,v) = g(j,j)$). This implies that $x^\frown i\leq e(u)$ and $x^\frown j \leq e(v)$, hence $\varepsilon(u,v) = (i,j)$. Finally
$$f(u,v) = g(\varepsilon(u,v)) = g(i,j) = \psi(g(i,i),g(j,j)) = \psi(f(u,u),f(v,v))$$
as desired.$\qed$\\

\section{Types}\label{sectiontypes}

In this section, we introduce the notion of type. We shall prove that each type corresponds to an equivalence class of analytic minimal strong gaps.

\begin{defn}
Let $n$ be a natural number. We call an $n$-type to the following collection of data:
\begin{enumerate}
\item A partition of $n$ into five sets, $n=A\cup B \cup C\cup D \cup E$
\item A function $\psi:\langle A \rangle^2\To B\cup\{\infty\}$ whose range covers $B$;
\item A partition $\mathcal{P}$ of $C$ into sets of cardinality either 1 or 2;
\item A function $\gamma:D \To B\cup E\cup\{\infty\}$ such that $|\gamma^{-1}(k)|\geq 2$ for every $k\in E$. 
\end{enumerate}
with the particularity that if $A=\emptyset$, then $B=D=E=\emptyset$ and all members of the partition $\mathcal{P}$ have cardinality 2.\\
\end{defn}

In principle, a type should be denoted by a tuple $(A,B,C,D,E,\psi,\mathcal{P},\gamma)$. However, we can denote the type by the shorter tuple $(\psi,\mathcal{P},\gamma)$ without mentioning explicitly $A$, $B$, $C$, $D$ and $E$, since this information is implicit. We point out some immediate consequences of the definition of type: the existence of $\psi$ implies that $B=\emptyset$ whenever $|A|\leq 1$. From the conditions on $\gamma$ it follows that $|D|\geq 2|E|$. And finally, if $A=\emptyset$ then $n$ must be even, since in that case we impose that $n=C$ and all members of $\mathcal{P}$ have cardinality 2.\\

Given an $n$-type $\alpha=(\psi,\mathcal{P},\gamma)$, we define an $n$-gap $\Gamma_{f^{\alpha}}$ as follows. We call $A^\ast = A$ if $A\neq\emptyset$, and $A^\ast = \{0\}$ if $A=\emptyset$. Let $M = A^\ast\cup \mathcal{P}\cup D$. We will have a function $f^\alpha:M^2\To n\cup\{\infty\}$. For $a\in \mathcal{P}$, let $\sigma(a) = \min(a)$ and $\tau(a) = \max(a)$ so that $a = \{\sigma(a),\tau(a)\}$ for every $a\in\mathcal{P}$. We define $\sigma(k) = k$ for $k\in A^\ast\cup D$, and $\tau(k) = \gamma(k)$ for $k\in D$. Notice that $\sigma:M\To n$ and $\tau:\mathcal{P}\cup D \To n\cup\{\infty\}$. The function $f = f^\alpha$ is defined as:
\begin{enumerate}
\item $f(i,i) = \sigma(i)$ for $i\in M$;
\item $f(i,j) = \psi(i,j)$ for $i,j\in A$, $i\neq j$;
\item $f(i,j) = \sigma(i)$ if $i\in M\setminus A^\ast$ and ($j\in A^\ast$ or $\sigma(i)<\sigma(j)$).
\item $f(i,j) = \tau(j)$ if $j\in M\setminus A^\ast$ and ($i\in A^\ast$ or $\sigma(i)>\sigma(j)$).\\
\end{enumerate}

Observe that $f$ satisfies the following properties:

\begin{itemize}
\item[(P1)] $A=pbranch(f)$;
\item[(P2)] $f$ is $\psi$-branch-reduced;
\item[(P3)] $k$ is $f$-attached to $l$ whenever $\{k,l\}\in\mathcal{P}$ (both when $k=l$ and when $k\neq l$);
\item[(P4)] $k$ is $f$-attached to $\gamma(k)$ for every $k\in D$.\\
\end{itemize}

\begin{thm}\label{thetypes}
For every $n$-type $\alpha$, $\Gamma_{f^\alpha}$ is a minimal strong $n$-gap. Moreover,
\begin{itemize}
\item Every minimal analytic strong $n$-gap is equivalent to $\Gamma_{f^\alpha}$ for some $n$-type $\alpha$;
\item If $\alpha$ and $\beta$ are different $n$-types, then $\Gamma_{f^\alpha}$ is not equivalent to $\Gamma_{f^\beta}$.
\end{itemize}
\end{thm}

Proof: We will combine a series of claims.\\

Claim 1: For every $n$-gap $\Gamma_g$ there exists some $n$-type $\alpha$ such that $\Gamma_{f^\alpha} \leq \Gamma_g$.\\

Proof of Claim 1.  We are given an $n$-gap $\Gamma_g$, with $g:m^2\To n\cup\{\infty\}$. Let $A = pbranch(g)$. For every $i\in A$, choose $u(i)\in m$ such that $g(u(i),u(i)) = i$. Define then $\psi:\langle A \rangle^2\To n\cup\{\infty\}$ as $\psi(i,j) = g(u(i),u(j))$, and then call $B$ to the range of the function $\psi$ excluding $\infty$. Let $F = n\setminus (A\cup B)$. We choose now a symmetric set $\Phi\subset\langle m\rangle^2$ (by symmetric we mean that $(i,j)\in \Phi$ implies $(j,i)\in \Phi$) such that $F\subset g(\Phi)$, but for every other symmetric $\Phi'\subset \Phi$, if $\Phi'\neq\Phi$ then $F\not\subset g(\Phi')$. We define a biniary relation $\approx$ on $n\cup\{\infty\}$ by saying that $k\approx l$ if there exists $(i,j)\in\Phi$ such that $g(i,j) = k$ and $g(j,i)=l$. Notice that this relation is symmetric ($k\approx l$ implies $l\approx k$) but not reflexive (we may have $k\not\approx k$), and also that for every $k\in F$ there exist at least one $l$ with $k\approx l$. Let now

$$ E = \{k\in F : \exists l,l', l\neq l'\in F  : k\approx l \text{ and } k\approx l'\},$$
$$ C = \{k\in F\setminus E : \exists l\in F\setminus E : k\approx l\}$$
$$ D = F\setminus (C\cup E).$$

Define the partition of $C$, $\mathcal{P} = \{\{k,l\} : k\approx l\}$, and given $k\in D$, define $\gamma(k)$ as the only element\footnote{Notice that, for $k\in D$, there is a unique choice for $\gamma(k)$. Namely, suppose that $k\approx l$ and $k\approx l'$ for $l\neq l'$. Since $k\not\in E$, either $l\not\in F$ or $l'\not\in F$. Say that $l\not\in F$, and $g(i,j)=k$, $g(j,i)=l$, $(i,j)\in\Phi$. Then $\Phi' = \Phi\setminus\{(i,j),(j,i)\}$ satisfies $F\subset g(\Phi')$ and contradicts the minimality in the choice of $\Phi$.} such that $k\approx \gamma(k)$. It is obvious that $\gamma(k) \not\in A = pbranch(g)$, and also $\gamma(k)\not\in F\setminus E$ since $k\not\in C$. Hence, we have $\gamma:D\To B\cup E\cup\{\infty\}$. So far, we collected all data in order to define our type $\alpha = (\psi,\mathcal{P},\gamma)$ and $f = f^\alpha$.\\

We notice that a special situation happens if $A=\emptyset$. In such case, by Lemma~\ref{pbranchdychotomy}, we can suppose without loss of generality that for every $k\in n$ there exists $l\in n\setminus\{k\}$ such that $l$ is $g$-attached to $k$. It is clear then that $B=\emptyset$, and using Lemma~\ref{attachproperties}, we get that $E=\emptyset$, that $C=F=n$, and that the partition $\mathcal{P}$ consists only of doubletons $\{k,l\}$, $k\neq l$, where $l$ is $g$-attached to $k$ and $k$ is $g$-attached to $l$.\\

 It remains to prove that $\Gamma_f\leq \Gamma_g$. We consider again $M = A^\ast\cup\mathcal{P}\cup D$, and the functions $\sigma$ and $\tau$ that we needed for the definition of $f = f^\alpha$. Let $q = |D\cup\mathcal{P}|+1$. We have to define an element $x\in m^{<q}$ and a function $e:M\To m^q$ so that the associated reduction map $\varepsilon:M^2\To m^2$ satisfies $f = g\circ \varepsilon$. Let us enumerate $D\cup\mathcal{P} = \{z_0,\ldots,z_{q-2}\}$ in such a way that $\sigma(z_0)<\sigma(z_1)<\cdots<\sigma(z_{q-2})$. For every $\xi < q-1$ we take $(i_\xi,j_\xi)\in \Phi$ such that $g(i_\xi,j_\xi) = \sigma(z_\xi)$, and moreover $g(j_\xi,i_\xi) = \gamma(z_\xi)$ when $z_\xi\in D$, while $\{g(j_\xi,i_\xi),\sigma(z_\xi)\}=z_\xi\in \mathcal{P}$ when $z_\xi\in \mathcal{P}$. In any case, $g(j_\xi,i_\xi) = \tau(z_\xi)$ for all $\xi$. Define then 
\begin{itemize}
\item $x= (j_0,j_1,\ldots,j_{q-2})$, \item $e(z_\xi) = (j_0,\ldots,j_{\xi-1},i_\xi,0,0,\ldots)\text{ for } \xi<q-2$, \item $e(k) = (j_0,j_1,\ldots,j_{q-2},u(k))\text{ for } k\in A^\ast.$\footnote{If $A=\emptyset$ and $A^\ast = \{0\}$, define $u(0)=0$}\\
 
\end{itemize}

\begin{center}
\includegraphics[scale=0.75]{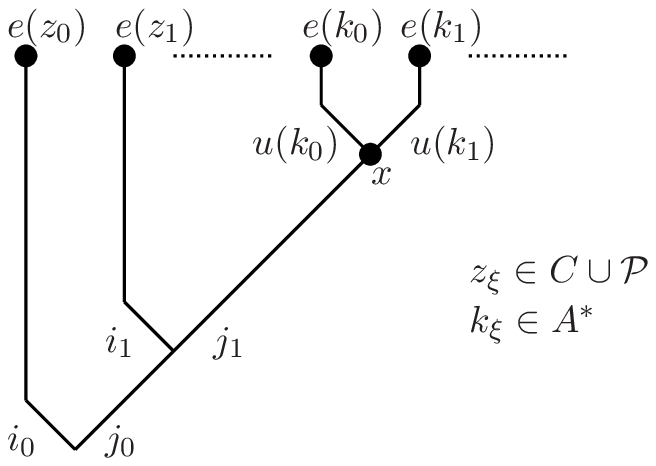}
\end{center}

It is straightforward to check that in this way, $f=g\circ \varepsilon$.\\

Claim 2: If $\Gamma_g$ is an $n$-gap with $\Gamma_g \leq \Gamma_{f^\alpha}$, then $pbranch(g) = pbranch(f^\alpha)$.\\

Proof of Claim 2: Let $f=f^\alpha$. By Lemma~\ref{hereditarypbranch}, we only have to prove that $pbranch(g)\subset A = pbranch(f)$. Let $g:m^2\To n\cup\{\infty\}$ and let $\varepsilon:m^2\To M^2$ be a reduction map coming from $x\in M^{<k}$ and $e:m\To M^k$. We distinguish two cases: $A\neq\emptyset$ and $A=\emptyset$.\\

If $A\neq\emptyset$, choose some $i\in A$. We pick $j\in pbranch(g)\setminus A$ and we work towards a contradiction. Since $j\in pbranch(g)$, there exists $p\in m$ such that $g(p,p) = j$. Remember that $$j=g(p,p) = f(\varepsilon(p,p)) = f(inc(e(p),x)).$$
Also, there must exist $q\in m$ such that $g(q,q) = i$ and again
$$i=g(q,q) = f(\varepsilon(q,q)) = f(inc(e(q),x)),$$
but in this case, since $i\in A$ we have that $f^{-1}(i) = \{(i,i)\}$, therefore $inc(e(q),x) = (i,i)$. That is: $x^\frown i\leq e(q)$. We distinguish two cases. First case: if $x\not\leq e(p)$, then $inc(e(p),x) = inc(e(p),e(q))$, therefore $j = f(inc(e(p),x)) = f(inc(e(p),e(q))) = g(p,q)$, which implies that $j\not\in pbranch(g)$ since $p\neq q$. This is a contradiction. Second case: $x\leq e(p)$. Call $(w,w)=inc(e(p),x)$. Notice that $w\neq i$ since $f(w,w) = j\not\in A$ while $f(i,i)=i\in A$. We have that $x^\frown w \leq e(p)$ and $x^\frown i \leq e(q)$, therefore $inc(e(p),e(q))= (w,i)$. Looking at item (1) of the definition of $f^\alpha=f$, we get that $f(w,w) = \sigma(w) = j\not\in A=A^\ast$. Loooking then at item (3) of the definition of $f$, we obtain that $f(w,i) = \sigma(w) = j$. Therefore $g(p,q) = f(inc(e(p),e(q))=f(w,i) = \sigma(w)=j$, which contradicts that $j\in pbranch(g)$.\\

We deal now with the case when $A=\emptyset$. Remember that in this case, we suppose that $n=C$ and all members of the partition $\mathcal{P}$ have cardinality 2. We first notice that for every $a\in\mathcal{P}$, $\sigma(a)\in pbranch(g)$ if and only if $\tau(a)\in pbranch(g)$. Namely, if $\sigma(a)\not\in pbranch(g)$, then $\sigma(a) = g(p,q) = f(\varepsilon(p,q))$ with $p\neq q$, and then $g(q,p) = \tau(a)$, hence $\tau(a)\not\in pbranch(g)$. Similarly if $\tau(a)\not\in pbranch(g)$, then $\sigma(a)\not\in pbranch(g)$. So suppose that $\sigma(a),\tau(a)\in pbranch(g)$, $\sigma(a) = g(p,p)$, $\tau(a) = g(q,q)$. Notice that $inc(e(p),x)\neq inc(e(q),x)$ because $f(inc(e(p),x)) = \sigma(a)$ and $f(inc(e(q),x)) = \tau(a)$. We distinguish three cases. \begin{itemize}

\item Case 1: $e(p)\wedge x \neq e(q)\wedge x$. Without loss of generality, we suppose that $e(p)\wedge x < e(q)\wedge x$. In this case, $g(p,q) = f(inc(e(p),e(q)) = f(inc(e(p),x)) =\sigma(a)$, which contradicts that $\sigma(a)\in pbranch(g)$.

\item Case 2: $e(p)\wedge x = e(q)\wedge x = x$. In this case, $inc(e(q),x)$ is a pair of the form $(w,w)$ and for all such pairs we have $f(w,w) = \sigma(w)$. This contradicts that $f(w,w) = f(inc(e(q),x)) = g(q,q)= \tau(a)$ while $\tau(a)\not\in\{\sigma(i) : i\in M\}$.

\item Case 3:  $y=e(p)\wedge x = e(q)\wedge x < x$. Suppose that $y^\frown i \leq e(p)$, $y^\frown j\leq e(q)$ and $y^\frown k \leq x$. Then $\tau(a) = g(q,q) = f(inc(e(q),x)) = f(j,k)$, so by the definition of $f$, we get $a=k$. Similarly, $f(i,k) = g(p,p) = \sigma(a)$ so again using the definition of $f$,we get $a=i$, a contradiction since $k\neq i$. This finishes the proof of Claim 2.\\  
\end{itemize}

Claim 3: If $\Gamma_{f^{\beta}} \leq \Gamma_{f^\alpha}$, then $\alpha=\beta$.\\

Proof of Claim 3: Let $\alpha = (\psi,\mathcal{P},\gamma)$ and $\beta = (\psi',\mathcal{P}',\gamma')$ with the obvious notations for the rest of objects associated with $\alpha$ and $\beta$. By Claim 2, $A = pbranch(f) = pbranch(f') = A'$. By Lemma~\ref{reducedpreserved}, $\psi = \psi'$, hence also $B=B'$. The case $A=\emptyset$ is quite simple because we only have to check that $\mathcal{P} = \mathcal{P}'$. This follows from Lemma~\ref{attachpreserved} since $\{k,l\}\in\mathcal{P}$ if and only if $k$ is $f$-attached to $l$ and vice-versa. For the case $A\neq\emptyset$,  we check the following items:\\

\begin{itemize}
\item $E\subset E'$. Pick $k\in E$. Since $|\gamma^{-1}(k)|\geq 2$, there exist $k_1,k_2\in D$, $k_1\neq k_2$, such that $k_1$ and $k_2$ are $f$-attached to $k$. By Lemma~\ref{attachpreserved}, $k_1$ and $k_2$ are also $f'$-attached to $k$. By Lemma~\ref{attachproperties}, this implies that $k$ cannot be $f'$-attached to any $l\in n\cup\{\infty\}$. From this, we get that $k\not\in C'\cup D'$. Also $k\not\in A'\cup B'=A\cup B$ since $k\in E$. We conclude that $k\in E'$.\\

\item $E'\subset E$. Pick $k\in E'$. Then, similarly as in the previous case, there exist  $k_1,k_2\in D'$, $k_1\neq k_2$, such that $k_1$ and $k_2$ are $f'$-attached to $k$, and this implies that $k$ is not $f'$-attached to any $l\in n\cup\{\infty\}$. Using Lemma~\ref{attachpreserved}, we get that $k$ is neither $f$-attached to any $l\in n\cup\{\infty\}$. This implies that $k\not\in C\cup D$. Also $k\not\in A\cup B = A'\cup B'$ since $k\in E'$. We conclude that $k\in E$.\\

\item $D\subset D'$. Pick $k\in D$. We only have to check that $k\not\in C'$. We know that $k$ is $f$-attached to $\gamma(k)\in B\cup E\cup\{\infty\}$. Hence, $k$ is also $f'$-attached to $\gamma(k)$. On the other hand, if $k\in C'$, then $k$ is $f'$-attached to some $l\in C'$.  If $k$ is $f'$-attached to both $\gamma(k)$ and $l$, we obtain that $\gamma(k)=l$. But this is a contradiction, because $\gamma(k)\in B'\cup E'\cup\{\infty\}$, while $l\in C'$.\\

\item $C\subset C'$. Pick $k\in C$. We only have to check that $k\not\in D'$. We know that there exists $l\in C$ such that $k$ is $f$-attached to $l$. Hence, $k$ is also $f'$-attached to $l$. On the other hand, if $k\in D'$, then $k$ is $f'$-attached to some $\gamma'(k)\in B'\cup E'\cup\{\infty\}$ . Again, we conclude that $l = \gamma'(k)$, but $l\in C$ while $\gamma'(k)\in B\cup E\cup\{\infty\}$, a contradiction.\\

\end{itemize}

We knew that  $A=A'$, $B=B'$ and $\psi=\psi'$ and from the four items above we obtain that $E=E'$, $C=C'$ and $D=D'$. For every $k$ in $D$, $\gamma(k)$ is the only $l$ such that $k$ is $f$-attached to $l$, so using Lemma~\ref{attachpreserved} we get that $\gamma=\gamma'$. Also, for $k,l\in C$, we have that $\{k,l\}\in\mathcal{P}$ if and only if $k$ is $f$-attached to $l$. So for the same reason, $\mathcal{P} = \mathcal{P}'$. This finishes the proof of Claim 3.\\

The theorem follows now straighforward from the claims. We prove that each $\Gamma_{f^\alpha}$ is minimal: Suppose that we have $\Gamma_g\leq \Gamma_{f^\alpha}$; Then by Claim 1, there exists $f^\beta$ such that $\Gamma_{f^\beta}\leq \Gamma_g$, and by Claim 3, $\alpha=\beta$. The fact that all analytic minimal strong $n$-gaps are equivalent to some $\Gamma_{f^\alpha}$ follows now from Claim 1 and Theorem~\ref{strongcontain}. Finally, by Claim 3, if $\alpha\neq \beta$, then the minimal strong $n$-gaps, $\Gamma_{f^\alpha}$ and $\Gamma_{f^\beta}$ are not equivalent.$\qed$\\

\section{Remarks}

We notice that although we have identified each of the classes of analytic minimal strong $n$-gaps, we did not provide a criterion if a given $n$-gap of the form $\Gamma_f$ is minimal equivalent to some $\Gamma_{f^\alpha}$, except just saying $\Gamma_f\leq \Gamma_{f^\alpha}$ and this is relatively satisfactory since we have a nice criterion by Lemma~\ref{criterion}. It is possible to give an intrinsic criterion but we found it too technical. The point is that properties $(P1)-(P4)$ enumerated after the definition of the functions $f^\alpha$ are necessary but not sufficient for $\Gamma_f$ to be minimal equivalent to $\Gamma_{f^\alpha}$. One thing is that for an $f$ satisfying those properties $(P1)-(P4)$ and $\Gamma_g\leq \Gamma_f$ we could have $pbranch(g)$ strictly larger than $pbranch(f)$. Even if this cannot happen, still we may not have a minimal strong gap, and we can even have phenomena like the following:\\

\begin{prop}
There exist $n$-gaps of the form $\Gamma_g$ which are not minimal strong $n$-gaps, but all minimal strong $n$-gaps below $\Gamma_g$ are equivalent.
\end{prop}

Proof: We consider the 4-gap $\Gamma_g$ given by $g:3^2\To 4$, where $g(i,i)=3$, and $g(i,j) = 3-i-j$ for $i\neq j$, $i,j=0,1,2$. Notice that
\begin{itemize}
\item $pbranch(g) = \{3\}$
\item $i$ is $g$-attached to $i$ for every $i=0,1,2$.\\
\end{itemize}

Claim: If $\Gamma_h\leq \Gamma_g$, then $pbranch(h) = \{3\}$. Proof of the claim: By Lemma~\ref{hereditarypbranch}, $3\in pbranch(h)$. Let us suppose that $i\in pbranch(h)$ for some $i<3$. Let $h:m^2\To 4$, $\varepsilon$ the reduction map with associated function $e$ and $x$. There is some $p$ such that $h(p,p) = g(inc(e(p),x)) = i$. Notice that if $g(u,v) = i$, then $u\neq v$, hence $x\not\leq e(p)$. There is another $q$ with $h(q,q) = g(inc(e(q),x))=3$ and now $g(u,v) = 3$ implies $u=v$, hence $x\leq e(q)$. Finally, $h(p,q) = g(inc(e(p),e(q))) = g(inc(e(p),x) = i$ and $p\neq q$, so $i\not\in pbranch(h)$.\\

Now, suppose $\Gamma_{f^\alpha}\leq \Gamma_h$, then $pbranch(f^\alpha) = \{3\}$ and $i$ is $f^\alpha$-attached to $i$ for every $i\in\{0,1,2\}$. This implies that the type $\alpha$ must be given by $A=\{3\}$, $C=\{0,1,2\}$, $\mathcal{P} = \{\{0\},\{1\},\{2\}\}$. This shows that there is only one possible type of minimal strong gap below $g$. Let us fix now $\alpha$ this type, and we suppose for contradiction that $\Gamma_g\leq \Gamma_{f^\alpha}$. For the type $\alpha$, $M = A\cup \mathcal{P} = \{3\}\cup\{\{0\},\{1\},\{2\}\}\equiv \{3,0,1,2\}$ and we can check that $f^\alpha:M^2\To 4$ is given by $f^\alpha(i,j) = \min(i,j)$. Let $\varepsilon$, $e:3\To 4^k$ and $x$ witnessing that $\Gamma_g\leq \Gamma_{f^\alpha}$. We distinguish two cases.

Case 1: $e(0)\wedge e(1) = e(0)\wedge e(2) = e(1)\wedge e(2)$. Say that $y^\frown {u_i} \leq e(i)$ for $i=0,1,2$. Then
$$g(i,j) = f^\alpha(e(i),e(j)) = f^\alpha(inc(e(i),e(j))) = f^\alpha(u_i,u_j) = \min(u_i,u_j),$$  
but this is not possible, because $g(i,j)$ takes values $0,1,2$ but $\min(u_i,u_j)$ takes only two values when $i=0,1,2$.

Case 2: If case 1 does not hold, then $inc(e(i_0),e(i_1)) = inc(e(i_0),e(i_2))$ for some reordering $\{0,1,2\} = \{i_0,i_1,i_2\}$. But then, since $f^\alpha$ is symmetric, again $g(i,j) = f^\alpha(inc(e(i),e(j)))$ takes only two possible values: $f^\alpha(inc(e(i_0),e(i_1)))$ and $f^\alpha(inc(e(i_1),e(i_2)))$, while $g(i,j)$ takes values $0,1,2$. Again, a contradiction. $\qed$\\

We point out that the key phenomenon is that, in general, for $f=f^\alpha$, if we take $W \subset M$ with at most one element in $A^\ast$, then there is an enumeration $W=\{w_0<w_1<\ldots\}$ such that $f(w_i,w_j) = f(w_i, w_k)$ and $f(w_j,w_i) = f(w_k,w_i)$ whenever $i<j<k$. This property is partially inherited to any $\Gamma_g\leq \Gamma_f$ under suitable hypotheses.

\section{The cases when $n=2,3$}

In this section we shall give the interpretation of our results for $2$-gaps and $3$-gaps. We will use a different notation for the $n$-gaps of the form $\Gamma_f$. Given $(i_1,j_1),\ldots,(i_k,j_k)\in m^2$, we will denote by $J^m_{(i_1j_1,\cdots, i_kj_k)}$ the ideal of $\mathcal{P}(m^{<\omega})$ generated by all $(i_\xi,j_\xi)$-combs, $\xi=1,\ldots,k$. In order to compute $f^\alpha:M\To n\cup\{\infty\}$ for a given type $\alpha$, we identify $M$ with a natural number $m$ by enumerating it.

\begin{thm}
There are six equivalence classes (four up to permutation), of analytic minimal strong 2-gaps.
\begin{enumerate}
\item $\{J^2_{(00,11,10)},J^2_{(01)}\}$,
\item $\{J^2_{(00)}, J^2_{(11)}\}$,
\item $\{J^2_{(00)}, J^2_{(11,01,10)}\}$ and its permutation,
\item $\{J^2_{(00)}, J^2_{(11,10)}\}$ and its permutation.
\end{enumerate}

\end{thm}

Proof: We just have to consider all possible 2-types $\alpha$ and construct the associated $\Gamma_{f^\alpha}$. Each case above corresponds with: 
\begin{enumerate}
\item $A=\emptyset$. We notice that $\{J^2_{(00,01)},J^2_{(11,10)}\}$ is in the same equivalence class.
\item $A=2$. Hence $B=\emptyset$, and $\psi$ is constant equal to $\infty$.
\item $|A|=1$ and $|C|=1$. The gap shown corresponds to $0\in A$ and $1\in C$, and conversely for its permutation. We notice that $\{J^2_{(00,11)},J^2_{(01,10)}\}$ is in the same equivalence class. 
\item $|A|=1$ and $|D|=1$. The gap shown corresponds to $0\in A$ and $1\in D$, and conversely for its permutation. Notice that $\gamma$ can only take the value $\infty$.$\qed$\\ 
\end{enumerate}

\begin{thm}
There are 31 equivalence classes of analytic minimal strong 3-gaps (nine up to permutation):
\begin{enumerate}
\item $\{J^3_{(00)}, J^3_{(11)}, J^3_{(22)}\}$.

\item $\{J^2_{(00)}, J^2_{(11)}, J^2_{(01,10)}\}$ and its three permutations.

\item $\{J^2_{(00)},J^2_{(11)}, J^2_{(01)}\}$ and its six permutations.

\item $\{J^3_{(00)}, J^3_{(11)}, J^3_{(22,02,20,21,12)}\}$ and its three permutations.

\item $\{J^3_{(00)}, J^3_{(11)}, J^3_{(22,20,21)}\}$ and its three permutations.

\item $\{J^2_{(00)}, J^2_{(11,10)}, J^2_{(01)}\}$ and its three permutations.

\item $\{J^3_{(00)}, J^3_{(11,01,10)}, J^3_{(22,02,20,12,21)}\}$ and its three permutations. 

\item $\{J^3_{(00)}, J^3_{(11,01,10)}, J^3_{(22,20,21)}\}$ and its six permutations.

\item $\{J^3_{(00)}, J^3_{(11,01)}, J^3_{(22,20,21)}\}$ and its three permuations.
\end{enumerate}

\end{thm}

Proof: We just have to consider all possible 3-types $\alpha$ and construct the associated $\Gamma_{f^\alpha}$. Each case above corresponds with:

\begin{enumerate}
\item $|A|=3$.
\item $|A|=2$, $|B|=1$ and $\psi$ is of the form $\psi(i,j) = \psi(j,i) \in B$. .
\item $|A|=2$, $|B|=1$ and $\psi$ is of the form $\psi(i,j)\in B$, $\psi(j,i)=\infty$.
\item $|A|=2$, $B=\emptyset$, $|C|=1$. Notice that $\psi$ must be constant equal to $\infty$.
\item  $|A|=2$, $B=\emptyset$, $|D|=1$. Notice that $\psi$ and $\gamma$  must be constant equal to $\infty$.
\item $|A| = 1$, $|C| = 2$, $|\mathcal{P}|=1$.
\item $|A| = 1$, $|C|=2$, $|\mathcal{P}|=2$.
\item $|A| = 1$, $|C|=1$, $|D|=1$. Notice that $\gamma$ must be constant equal to $\infty$.
\item $|A| = 1$, $|D|=2$. Notice that $\gamma$ must be constant equal to $\infty$.$\qed$\\
\end{enumerate}

\section{An application}

In our previous work \cite{multiplegaps} we defined a multiple gap $\Gamma$ to be a clover if it is not possible to find a set $a$, a subset $B\subset n$ and $i\not\in B$ such that $a\in \Gamma(i)^\perp$, but $\{\Gamma(j)|_a : j\in B\}$ is a multiple gap. We constructed one clover but we left as an open a question whether analytic clovers exist or not. Our results allow to answer this question in the case of strong multiple gaps, since we only have to check the minimal ones.

\begin{prop}
Suppose that $\Gamma$ is an analytic strong $n$-gap with $n\geq 3$. Then there is a partition $n=X\cup Y$ and $a\subset \omega$ such that
\begin{enumerate}
\item $a\in \Gamma(i)^\perp$ for $i\in Y$
\item $|X|\geq 2$ and $\{\Gamma(j)|_a : j\in X\}$ is a strong multiple gap.
\end{enumerate}
In particular, $\Gamma$ cannot be a clover.
\end{prop}

Proof: It is enough to consider the case when $\Gamma = \Gamma_{f^\alpha}$ is a minimal strong gap. There are many ways of producing such a decomposition. For instance, if we consider $a =A^{<\omega}\subset M^{<\omega}$, then this $a$ has that property for $X=A\cup B$ and $Y = C\cup D\cup E$. Another way: if we pick $i\in A$, and $a$ consists of all tuples $(s_0,\ldots,s_{2k-1})$ such that $s_{2p}\neq i$ for every $p<k$, then $a$ is a subtree\footnote{The set $a$ can be identified with the tree $N^{<\omega}$, where $N = (n\setminus\{i\})\times n$.} of $M^{<\omega}$ where all possible combs and chains appear except $i$-chains. Hence it satifies the property above for $X=n\setminus\{i\}$ and $Y=\{i\}$. We can do that whenever $A\neq\emptyset$. If $A=\emptyset$ and we choose $P=\{u,v\}\in\mathcal{P}$, then $a=\{0,P\}^{<\omega}$ also works for $X=\{u,v,0\}$ and $Y=n\setminus\{u,v,0\}$. $\qed$

%
%
%

\end{document}